\documentclass[11pt]{article}
\usepackage{amsmath,amssymb,amscd,epsfig,graphics}
\def\mapdownright#1{\Big\downarrow\rlap{$\vcenter{\hbox{$\scriptstyle#1$}}$}}

\newtheorem{theorem}{Theorem}

\newtheorem{Claim}{Claim}[section]

\newcommand{\prf}{\vspace{.05in}
                    \noindent {\sc Proof} \hspace{.05in}}

\newcommand{\bp}{{\mathbb P}}

\newcommand{\bx}{{\mathcal X}}
\newcommand{\bl}{{\mathcal L}} 
 
\newcommand{\bo}{{\mathcal O}}

\newcommand{\restr}[1]{\!\!\mid_{\,#1}}

\newcommand{\im}{\mathop{\rm Im}\nolimits}

\renewcommand{\Box}{\square}

\begin{document}

\begin{center}
{\LARGE Linear Systems of Plane Curves with a Composite Number of Base Points of Equal Multiplicity} \\
\vspace{0.2in}
{\large Anita Buckley and Marina Zompatori}\\
\vspace{0.15in}
{\large 28. June 2001}
\vspace{0.1in}
\end{center}

\setcounter{footnote}{0}

\renewcommand{\thefootnote}{\arabic{footnote}}

{\small
\begin{center} {\sc abstract} \end{center}
{\leftskip=30pt \rightskip=30pt
In this article we study  linear systems of plane curves of degree $d$ passing through general base points with the same mutliplicity at each of them. These systems are known as \textit{homogeneous} linear systems. We especially investigate for which of these systems, the base points, with their multiplicities, impose independent conditions and which homogeneous systems are empty. Such systems are called \textit{non--special}.  We extend the range of homogeneous linear  systems that are known to be non--special.
A theorem of Evain states that the systems of curves of degree $d$ with $4^h$ base points with equal multiplicity are non--special.  The analogous result for $9^h$ points was conjectured. Both of these will follow, as corollaries, from the main theorem proved in this paper. Also, the case of $4^{h}9^{k}$ points will follow from our result.
The proof uses a degeneration technique developed by C. Ciliberto and R. Miranda.
\footnote{We wish to thank the organizers of Pragmatic 2001 for sponsoring  our stay in Catania, Rick Miranda and Ciro Ciliberto for introducing us to problems on linear systems and for many invaluable conversations. We would also like to thank Dan Abramovich and Bal\'{a}zs Szendr\H{o}i  for corrections and helpful comments.}.  \par}}

\section*{Introduction}
\label{introduction}

Fix the projective plane $\bp^{2}$ and $n$ general points $q_1,\ldots,q_n$ on it. Consider the linear system consisting of plane curves of degree $d$ with multiplicity at least $m_i$ at each point $q_i$. In this paper we will restrict ourselves 
to linear systems of curves which pass through every point with the same fixed multiplicity $m,$ i.e. 
$m_i = m$ for all $i = 1,\ldots, n.$ Such systems are called \textit{homogeneous} and are denoted by  $\bl_d(m^n).$

\pagebreak 
The linear system of all plane curves of degree $d$ has projective dimension 
$$\left(\begin{array}{c}
d+2\\
2        
\end{array}
\right)-1=\frac{d(d+3)}{2}.$$
Each base point of multiplicity $m$ imposes $\frac{m(m+1)}{2}$ conditions. These conditions come from the vanishing of the coefficients in the Taylor expansion of the equation of the curve. 
  Define the \textit{virtual dimension} of $\bl_d(m^n)$ by
$$v=v(d,m,n)= \frac{d(d+3)}{2}-n\ \frac{m(m+1)}{2}.$$
The actual dimension of this system can not be less than $-1$ and if the dimension is equal to $-1,$ then the system is empty. Therefore  we define the \textit{expected dimension} to be
$$e=e(d,m,n)=\max\{-1,\ v\}.$$

Since we chose the points  $q_1,\ldots,q_n$ in general position, the dimension of
$\bl_d(m^n)$
achieves its minimum value. We call this the \textit{dimension} of $\bl_d(m^n)$ and denote it by $l.$
Observe that
$l\geq e.$  Equality implies that either the system is empty, or the conditions imposed by the multiple points are independent. 

The system $\bl_d(m^n)$
is \textit{non--special} if $l=e$. Otherwise, we say that the system is \textit{special.}\\

The main result of this paper is the following theorem.

\begin{theorem}
\label{thm}
If all the homogeneous linear systems, of all degrees and multiplicities, through $n_1$ points in general position are non--special and all the homogeneous systems through $n_2$ points in general position are non--special, then the linear systems $\bl_d(m^{n_1n_2})$ are  also non--special, for all $d$ and $m.$
\end{theorem}

We are actually able to prove more. The following is a slightly more precise version of the previous theorem, which actually uses the full strength of our proof.

\begin{theorem}
\label{thm2}
For fixed $d, m, n_{1}, n_{2}$ as above, and for $k$ an integer between

$$\frac{1}{2}\left(-3 + \sqrt{1 + 4n_{2}m(m + 1)}\right) \,\,\,\hbox{and}\,\,\,\frac{1}{2}\left(-1 + \sqrt{1 + 4n_{2}m(m + 1)}\right), $$
the linear system $\bl_d(m^{n_1n_2})$ will be non-special if systems
$$\bl_d((k + 1)^{n_{1}}), \bl_d(k^{n_{1}}), \bl_k-1(m^{n_2}) \,\,\hbox{and}\,\, \bl_k(m^{n_2})$$ are non-special.
\end{theorem}

This is a generalization of Evain's theorem \cite{evain} in which 
it is proved that all systems of plane curves of degree $d$ through $4^h$ points with homogeneous multiplicity $m$ are 
non--special.
We prove our result by using a degeneration technique developed by C. Ciliberto and R. Miranda\cite{cil_mir} which we now describe.
 
\section{Degeneration Technique}
\label{degen}

We will now describe the degeneration technique needed for the proof. We first consider a degeneration of the plane $\bp^{2},$  then a degeneration of the bundle and finally of the base points.\\

Let us  start by degenerating $\bp^{2}.$ Let $\Delta$ be a complex disc around the origin. The product $V = \bp^{2} \times \Delta$ is equipped with two projections $p_{1} : V \rightarrow \Delta$ and $p_{2} : V \rightarrow \bp^{2}.$ We denote $V_{t} = \bp^{2} \times \{t\}.$

Consider $n_{1}$ general points in the plane $V_{0}$ and blow up $V$ at these points. We get a new three-fold $X$ with maps $f: X \rightarrow V,  \pi_{1} = p_{1} \circ f: X \rightarrow \Delta$ and $\pi_{2} = p_{2} \circ f: X \rightarrow \bp^{2}.$
The map $\pi_{1}$ gives a flat family of surfaces over $\Delta.$

$$ \begin{CD}
X                            \\
@VV{f}V                         \\
V   @>p_{2}>> \bp^{2}         \\
@VV{p_{1}}V                         \\
{\Delta}
\end{CD}$$

Let $X_{t}$ be the fiber of $\pi_{1}$ over $t \in \Delta.$ If $t\not= 0,$ then $X_{t}\cong V_{t}$ is a plane $\bp^{2}.$
By contrast, $X_{0}$ is the union of the proper transform of $V_{0},$ which we denote by $Y,$ and of $n_{1}$ exceptional divisors $\mathbf{P}_{i}$ for $i = 1, \ldots, n_{1}.$
Each $\mathbf{P}_{i}$ is isomorphic to the plane $\bp^{2}.$ The variety $Y$ is $\bp^{2}$ blown up at $n_{1}$ points. Denote this blow up by $b: Y\rightarrow \bp^2.$
Each $\mathbf{P}_{i}$ intersects $Y$ transversally along a curve $E_{i},$ which is a line in $\mathbf{P}_{i}$ and an exceptional divisor on $Y.$\\

Next we consider a degeneration of the bundle $\bo_{\bp^2}(d)$. 
Define a line bundle on $X$ by
$$\bo_X (d, k) =\pi_{2}^{*} \bo_{\bp^{2}} (d) \otimes \bo_X (kY).$$
The restriction of $\bo_X (d, k)$ to $X_{t},$ for $t \not= 0,$ is isomorphic to $\bo_{\bp^{2}}(d).$

Let $\bx (d, k)$ be restriction of $\bo_X (d, k)$ to $X_{0}.$ 
$\bx (d, k)$ is a flat limit of the bundle $\bo_{\bp^{2}}(d)$ on the general fiber $X_{t}$.
On each $\mathbf{P}_{i},$ the bundle $\bx (d, k)$ equals to $\bo_X(d,k)\restr{\mathbf{P}_{i}}$ which is isomorphic to $ \bo_{\bp^2}(k)$, since $Y$ intersects $\mathbf{P}_{i}$  along a line. 
On $Y$ bundle $\bx(d,k)$ is tensor product of

$$\pi_{2}^{*} \bo_{\bp^{2}}(d)\restr{Y}\cong \bo_Y(b^*(dL))$$  
 
and of 
$$\bo_X(kY)\restr{Y}\cong \bo_X(kY-kX_t)\restr{Y}\cong\bo_X(-\sum_{i=1}^{n_1} k\mathbf{P}_i)\restr{Y}\cong 
\bo_X(-\sum_{i=1}^{n_1} kE_i). $$
The last statement holds since $X_t\sim Y+\sum_1^{n_1}\mathbf{P}_i$ as divisors on $X$, and $\mathbf{P}_{i}$ intersects $Y$
along exceptional line.
So $$\bx(d,k)\restr{Y}\cong\bo_Y \left( dL - \sum_{i = 1}^{n_{1}} k E_{i} \right)$$ where we abuse notation and write $L$ for $b^{*}L$ on $Y.$\\

Now fix a positive integer  $n_{2}$ and put $n_{2}$ general points 
$p_{i,1},\ldots, p_{i,n_2}$ in each $\mathbf{P}_{i}.$
We can view these points as limits of $n_1n_2$ general points $p_{1,t},\ldots,p_{n_1n_2,t}$ in $X_t,\ t\neq 0.$
Consider then the linear system $\bl_t$ which is the system $\bl_d(m^{n_1n_2})$ in $X_t\cong\bp^2$ 
based at the points $p_{1,t},\ldots,p_{n_1n_2,t}$.

We also consider the system $\bl_0(d, k, m, n_1n_2)$ on $X_0$ consisting of divisors in $|\bx (d, k)|$ which vanish with  multiplicity at least $m$ at $p_{1,1},\ldots,p_{1,n_2},\ldots, p_{n_1,n_2}$.
As discussed above, $\bl_0 (d, k, m, n_{1}n_{2})$ on $X_0$  can be seen as a flat limit of the system 
$\bl_t=\bl_d(m^{n_1n_2})$.\\

By upper-semicontinuity, the dimension of $\bl_0 (d, k, m, n_{1}n_{2})$ is at least equal to the dimension of $\bl_d (m^{n_{1}n_{2}}).$
Therefore in order  to prove that the dimension of $\bl_{d}(m^{n_{1}n_{2}})$ is equal to the expected dimension, it is enough to prove  that, under the assumptions of Theorem~\ref{thm2}, 
$$\dim \bl_0 (d, k, m, n_{1}n_{2})=e(d, m, n_1n_2).$$

\section{The Transversality of the Restricted Systems}
\label{trans}

Recall from Section~\ref{degen} the line bundle $\bx(d,k)$ on $X_0$ such that  
$$\bx(d,k)\restr{Y}\cong \bo_Y\left( dL -\sum_{i=1}^{n_1}kE_i\right) $$ and
$$\bx(d,k)\restr{\mathbf{P}_{i}}\cong \bo_{\bp^2}(kL),$$ 
and that we put $n_2$ points with the same multiplicity $m$ in each $\mathbf{P}_i$. 
Global sections of $\bx(d,k)\restr{Y}$ correspond to curves in $\bp^2$ which pass through each of the $n_1$ points  with multiplicity at least $k$.
This implies that the restriction of the system $\bl_0(d,k,m,n_1n_2)$ to $Y$ is a system $\bl_Y$ isomorphic to $\bl_d(k^{n_1})$. In the same way,
the restriction of $\bl_0(d,k,m,n_1n_2)$ to $\mathbf{P}_{i}\cong\bp^2$ is a system $\bl_{\mathbf{P}_i}$ 
isomorphic to $\bl_k(m^{n_2})$.

The system $\bl_0$ is a linear system on a reducible scheme $X_0$. Its elements consist of $\alpha\in\bl_Y$ on $Y$ and of $\beta_i\in\bl_{\mathbf{P}_i}$ on each $\mathbf{P}_{i}$ which restrict to the same divisor on the lines $E_{i}.$ 

In order to compute the dimensions of linear systems, it is easier to use the dimensions of the corresponding vector spaces.
We denote by $\mathbf{L}_{d}(k^{n_{1}})$ the vector space whose projectivization is $\bl_{d}(k^{n_{1}})$ and similarly we introduce the vector spaces $\mathbf{L}_{k}(m^{n_{2}}), \mathbf{L}_{Y}, \mathbf{L}_{\mathbf{P}_{i}} \hbox{and} \mathbf{L}_{0}$ such that

$$\begin{array}{l}
\bp(\mathbf{L}_{d}(k^{n_{1}}))=\bl_{d}(k^{n_{1}}),\\
\bp(\mathbf{L}_{k}(m^{n_{2}}))=\bl_{k}(m^{n_{2}}),\\
\bp(\mathbf{L}_Y)=\bl_Y,\\
\bp(\mathbf{L}_{\mathbf{P}_i})=\bl_{\mathbf{P}_i},\\
\bp(\mathbf{L}_0)=\bl_0(d,k,m,n_1n_2).
\end{array}$$

By restriction to the lines $\cup_{i=1}^{n_1}E_i$ we have maps
$$\mathbf{L}_Y\stackrel{\rho_Y}{\longrightarrow}\oplus_{i=1}^{n_1}H^{0}(E_i,\bo_{E_i}(k))
\mbox{ and }
\mathbf{L}_{\mathbf{P}_i}\stackrel{r_i}{\longrightarrow}H^{0}(E_i,\bo_{E_i}(k)).$$
At the level of vector spaces,
$\mathbf{L}_0$ is the fibered product of $\mathbf{L}_Y$ and $\oplus_{i=1}^{n_1}\mathbf{L}_{\mathbf{P}_i}$ over
the vector space of restricted system on $\cup_{i=1}^{n_1}E_i$, which is 
$$\oplus_{i=1}^{n_1}H^{0}(E_i,\bo_{E_i}(k)).$$
This means
$$\mathbf{L}_0=\mathbf{L}_Y\times_{\oplus_{i=1}^{n_1}H^{0}(E_i,\bo (k))}\left( \oplus_{i=1}^{n_1}\mathbf{L}_{\mathbf{P}_i}\right) .$$

The situation is shown in the diagram below:
$$\begin{array}{rcccccccc}
\mathbf{L}_0 & \longrightarrow & \mathbf{L}_{\mathbf{P}_1} & \oplus & \mathbf{L}_{\mathbf{P}_2} & \oplus & \cdots & \oplus & 
\mathbf{L}_{\mathbf{P}_{n_1}} \\
\downarrow & & \mapdownright{r_1} & & \mapdownright{r_2}& & & & \mapdownright{r_{n_1}} \\
\mathbf{L}_Y &\stackrel{\rho_Y}{\longrightarrow} & H^{0}(E_1,\bo_{E_1}(k)) & \oplus & H^{0}(E_2,\bo_{E_2}(k)) & \oplus & 
\cdots & \oplus & H^{0}(E_{n_1},\bo_{E_{n_1}}(k)). 
\end{array}$$

Our task is to calculate the dimension of $\mathbf{L}_0$. 
By the definition of the fibred product
$$\mathbf{L}_0=\left\{ (\alpha,\beta)\in \mathbf{L}_Y\times \left( \oplus_{i=1}^{n_1}\mathbf{L}_{\mathbf{P}_i}\right) \ :\ 
\alpha\restr{\oplus_{i=1}^{n_1}E_i}=\beta \restr{\oplus_{i=1}^{n_1}E_i}\right\} .$$
An element of $\mathbf{L}_0 $ is determined by first choosing an element $\gamma$ in
$\im(\rho_Y)\cap\im(r_1,\ldots,r_{n_1})$
and then taking its inverse images $\alpha\in\mathbf{L}_Y$ and $\beta\in\left(\oplus_{i=1}^{n_1}\mathbf{L}_{\mathbf{P}_i}\right).$ 
Once we fix $\gamma$, the choice of $\alpha$ depends on 
$\dim\ker(\rho_Y)$ parameters and similarly the choice of $\beta$ depends on 
$\dim\ker(r_1,\ldots,r_{n_1})$ parameters.\\
\noindent From the above considerations we get
$$ \begin{array}{l}\dim\mathbf{L}_0=
\dim(\im(\rho_Y)\cap\im(r_1,\ldots,r_{n_1}))+\dim\ker(\rho_Y)+\dim\ker(r_1,\ldots,r_{n_1})\,\,.(*)
\end{array}$$

As  proved in~\cite{a_m}, we can apply the Generalized Transversality Lemma to the above diagram of vector spaces. By this lemma we can assume that
$\im(\rho_Y)$ and $\im(r_1,\ldots,r_{n_1})$ meet properly.
This means that
$$\begin{array}{l} 
\dim(\im(\rho_Y)\cap\im(r_1,\ldots,r_{n_1}))=\\
\max\left\{\dim\im(\rho_Y)+\dim\im(r_1,\ldots,r_{n_1})-n_1(k+1),\ 0\right\}. 
\end{array}$$

We can describe divisors in the kernel of $\rho_Y$ as elements in $\mathbf{L}_{Y}$ 
which contain  $\cup_{i=1}^{n_1}E_i$ as a component. These correspond to curvesof degree $d$ in $\bp^2$ which pass through each point of the blow up with multiplicity at least $k+1$.
This means $\ker(\rho_Y)\cong\mathbf{L}_d((k+1)^{n_1}).$
In the same way, divisors in the kernel of $r_i$ are elements of $\mathbf{L}_{\mathbf{P}_{i}}$ that contain $E_i$ as component. These correspond to curves of degree $k-1$ in $\mathbf{P}_i$ and therefore $\ker(r_i)\cong \mathbf{L}_{k-1}(m^{n_2}).$\\
This gives rise to left exact sequences 
$$ \begin{array}{c}
0\longrightarrow \mathbf{L}_d((k+1)^{n_1})\longrightarrow \mathbf{L}_Y\stackrel{\rho_Y}{\longrightarrow} \oplus_{i=1}^{n_1}H^{0}(\bo_{E_i}(k)) \\
 \\
0\longrightarrow \mathbf{L}_{k-1}(m^{n_2})\longrightarrow \mathbf{L}_{\mathbf{P}_i}\stackrel{r_i}{\longrightarrow} H^{0}(\bo_{E_i}(k)).
\end{array}$$
Recall $\bl_Y\cong\bl_d(k^{n_1})$ and $\bl_{\mathbf{P}_i}\cong\bl_k(m^{n_2})$. From this we can calculate
$$\begin{array}{l}
\dim\im(\rho_Y)=\dim\mathbf{L}_d(k^{n_1})-\dim\mathbf{L}_d((k+1)^{n_1}),\\
\dim\im(r_1,\ldots,r_{n_1})=n_1(\dim\mathbf{L}_k(m^{n_2})-\dim\mathbf{L}_{k-1}(m^{n_2}))\\
\dim\ker(\rho_Y)=\dim\mathbf{L}_d((k+1)^{n_1}),\\
\dim\ker(r_1,\ldots,r_{n_1})=n_1\dim\mathbf{L}_{k-1}(m^{n_2}).
\end{array}$$

Since, by assumption in our theorem, we know that the dimension of all homogeneous systems through $n_1$ or $n_2$ points is equal to the expected dimension, we are able to compute the dimensions of all the kernels and images above, and therefore we are able to compute the dimension of $\mathbf{L}_{0}.$

\section{Dimension of the Fibered Product $\mathbf{L}_0$}
\label{comp}

In the sequel we assume that  $n_1$ and $n_2$ are such that linear systems $\bl_d(m^{n_1})$ and $\bl_d(m^{n_2})$ are non-special for every $d$ and $m$. By definition of non-special system, the corresponding vector spaces have dimension
$$\dim\mathbf{L}_d(m^{n_1})=\max\left\{\frac{d(d+3)}{2}+1-n_1\frac{m(m+1)}{2},\ 0\right\}$$
and
$$\dim\mathbf{L}_d(m^{n_2})=\max\left\{\frac{d(d+3)}{2}+1-n_2\frac{m(m+1)}{2},\ 0\right\}.$$
Our task is to prove that then the linear system $\bl_d(m^{n_1n_2})$ is also non-special. As mentioned in Section~\ref{degen}, the dimension of the vector space $\mathbf{L}_0$ on $X_0$ by semi-continuity satisfies the following inequality
$$\dim(\mathbf{L}_0)\geq\dim\mathbf{L}_d(m^{n_1n_2}).$$
Therefore it is enough to prove
$$\dim\mathbf{L}_0=\max\left\{\frac{d(d+3)}{2}+1-n_1n_2\frac{m(m+1)}{2},\ 0\right\}.$$\\

Now we are ready to start the computation.
\begin{Claim}
In the computation for the dimension of $\mathbf{L}_{0},$ given $d, k, m, n_1, n_2$ as above, we have:

$$\begin{array}{l}
\dim\ker(\rho_Y)+\dim\ker(r_1,\ldots,r_{n_1})=\\
 \\
\dim\mathbf{L}_d((k+1)^{n_1})+n_1\dim\mathbf{L}_{k-1}(m^{n_2})=
\end{array}$$
$$\left\{\begin{array}{ll}
\frac{d(d+3)}{2}+1-n_1n_2\frac{m(m+1)}{2}-n_1(k+1) & \hbox{ if } \dim\mathbf{L}_d((k+1)^{n_1})\geq 0,\\ 
 & \hbox{ }\,\,\,\,\,\dim\mathbf{L}_{k-1}(m^{n_2})\geq 0\\
 & \\
\frac{d(d+3)}{2}+1-n_1\frac{(k+1)(k+2)}{2} & \hbox{ if } \dim\mathbf{L}_d((k+1)^{n_1})\geq 0,\\ 
& \hbox{ }\,\,\,\,\, \dim\mathbf{L}_{k-1}(m^{n_2})=0\\
 & \\
n_1\left(\frac{(k-1)(k+2)}{2}+1-n_2\frac{m(m+1)}{2}\right)& \hbox{ if } \dim\mathbf{L}_d((k+1)^{n_1})=0,\\  
& \hbox{ }\,\,\,\,\, \dim\mathbf{L}_{k-1}(m^{n_2})\geq 0\\
 & \\
0 & \hbox{ if } \dim\mathbf{L}_d((k+1)^{n_1})=0,\\
& \hbox{ } \,\,\,\,\, \dim\mathbf{L}_{k-1}(m^{n_2})=0.
\end{array}\right. $$\\
\end{Claim}
\prf
Suppose $ \dim\mathbf{L}_d((k+1)^{n_1})\geq 0 $ and $\dim\mathbf{L}_{k-1}(m^{n_2})\geq 0,$ then 
$\dim \mathbf{L}_{d}((k + 1)^{n_{1}}) + n_{1} \dim \mathbf{L}_{k - 1}(m^{n_{2}})$ is equal to
$$ \frac{d(d + 3)}{2} + 1 - n_{1}\frac{(k + 1)(k + 2)}{2} + n_{1}\frac{(k - 1)(k + 2)}{2} + n_{1} - n_{1}n_{2}\frac{m(m + 1)}{2}$$
which simplifies to the given expression.
This proves the first case. The other cases follow easily, given the definition of expected dimension, since one or both of the terms will be zero. \hfill $\Box$\\

Before we continue observe 
$$\dim\mathbf{L}_d(k^{n_1})\geq\dim\mathbf{L}_d((k+1)^{n_1})\hbox{ and } \dim\mathbf{L}_k(m^{n_2})\geq\dim\mathbf{L}_{k-1}(m^{n_2}),$$
and  $\dim\mathbf{L}_d(m^{n_{i}})=0 \hbox{ if and only if } \frac{d(d+3)}{2}+1-n_{i}\frac{m(m+1)}{2}\leq 0$ for $i = 1, 2.$
Bearing this in mind we get
$$\begin{array}{l}
\dim(\im(\rho_Y)\cap\im(r_1,\ldots,r_{n_1}))=\\
 \\
= \max\left\{\dim\im(\rho_Y)+\dim\im(r_1,\ldots,r_{n_1})-n_1(k+1),\ 0\right\}\\
 \\
= \max\left\{\dim\mathbf{L}_d(k^{n_1})-\dim\mathbf{L}_d((k+1)^{n_1})+n_1(\dim\mathbf{L}_k(m^{n_2})-\dim\mathbf{L}_{k-1}(m^{n_2})) - n_{1}(k + 1),\ 0\right\}.
\end{array}$$
We can now formulate the following claim.
\begin{Claim}
The dimension of $\dim(\im(\rho_Y)\cap\im(r_1,\ldots,r_{n_1}))$ is equal to
$$\left\{\begin{array}{ll}
n_1(k+1) & \hbox{if } \dim\mathbf{L}_d((k+1)^{n_1})\geq 0,\\
& \hbox{ }\,\,\,\,\, \dim\mathbf{L}_{k-1}(m^{n_2})\geq 0;\\
 & \\
n_1\left(\frac{k(k+3)}{2}+1-n_2\frac{m(m+1)}{2}\right)& \hbox {if }\dim\mathbf{L}_d((k+1)^{n_1})\geq 0,\\
 & \hbox{ }\,\,\,\,\,\dim\mathbf{L}_k(m^{n_2})\geq 0,\ \dim\mathbf{L}_{k-1}(m^{n_2})=0; \\
 & \\
\frac{d(d+3)}{2}+1-n_1\frac{k(k+1)}{2} & \hbox{if }\dim\mathbf{L}_d(k^{n_1})\geq 0,\\
& \hbox{ }\,\,\,\,\, \dim\mathbf{L}_d((k+1)^{n_1})=0, \dim\mathbf{L}_{k-1}(m^{n_2})\geq 0;\\  
 & \\
\frac{d(d+3)}{2}+1-n_1n_2\frac{m(m+1)}{2}& \hbox{if }\dim\mathbf{L}_d((k+1)^{n_1})=0,\\
& \hbox{ }\,\,\,\,\, \dim\mathbf{L}_{k-1}(m^{n_2})=0, \frac{d(d+3)}{2}+1-n_1n_2\frac{m(m+1)}{2}\geq 0;\\
 & \\
0 & \hbox { otherwise. }
\end{array}\right. $$\\
\end{Claim}
\prf
First consider the case $\dim \mathbf{L}_{d}((k + 1)^{n_{1}})\geq 0$ and $\dim \mathbf{L}_{k - 1}(m^{n_{2}}) \geq 0.$
Under this hypothesis 
$$\dim\mathbf{L}_d(k^{n_1})-\dim\mathbf{L}_d((k+1)^{n_1})+n_1(\dim\mathbf{L}_k(m^{n_2})-\dim\mathbf{L}_{k-1}(m^{n_2})) - n_{1}(k + 1)$$ is equal to
$\frac{d(d + 3)}{2} + 1 - n_{1}\frac{k(k + 1)}{2} - \frac{d(d + 3)}{2} - 1 + n_{1}\frac{k(k + 3)}{2} + n_{1} - n_{1}n_{2}\frac{m(m + 1)}{2} - n_{1}\frac{(k - 1)(k + 2)}{2} - n_{1} + n_{1}n_{2}\frac{m(m + 1)}{2} - n_{1}(k + 1).$
All the terms containing $d$ disappear and we are left with
$$\frac{n_{1}}{2}(- k^{2} - k + k^{2} + 3k + 2 + k^{2} + 3k - k^{2} - k + 2 - 2k - 2)$$ 
which is exactly $n_{1}(k + 1).$ In particular it is non-negative.

In the second case $\dim \mathbf{L}_{k - 1}(m^{n_{2}}) = 0,$ the expression 
$$\dim \mathbf{L}_{d}(k^{n_{1}}) -  \dim \mathbf{L}_{d}((k + 1)^{n_{1}}) + n_{1} \dim\mathbf{L}_{k}(m^{n_{2}}) - n_{1}(k + 1)$$ gives 
$ - n_{1}\frac{k(k + 1)}{2} + n_{1}\frac{(k + 1)(k + 2)}{2} + n_{1}\left(\frac{k(k + 3)}{2} + 1 - 
n_{2}\frac{m(m + 1)}{2}\right) - n_{1}(k + 1),$
which is exactly equal to $n_{1}\left(\frac{k(k + 3)}{2} + 1 - n{2}\frac{m(m + 1)}{2}\right).$
By hypothesis, this quantity is non-negative.

In the third case, we assume that $\dim\mathbf{L}_{d}(k^{n_{1}})$ and $\dim\mathbf{L}_{d - 1}(m^{n_{2}})$ are both non-negative and that $\dim\mathbf{L}_{d}((k + 1)^{n_{1}}) = 0.$
So the expression 
$$\dim\mathbf{L}_d(k^{n_1})-\dim\mathbf{L}_d((k+1)^{n_1})+n_1(\dim\mathbf{L}_k(m^{n_2})-\dim\mathbf{L}_{k-1}(m^{n_2})) - n_{1}(k + 1)$$ 
under these hypotheses gives

$$\begin{array}{l}
\frac{d(d + 3)}{2} + 1 - n_{1}\frac{k(k + 1)}{2} \\
+ n_{1}\left(\frac{k(k + 3)}{2} + 1 -n_{2}\frac{m(m + 1)}{2} - 
\frac{(k - 1)(k + 2)}{2} - 1 + n_{2}\frac{m(m + 1)}{2}\right) - n_{1}(k + 1),
\end{array}$$ 
which is equal to $\frac{d(d + 3)}{2} + 1 + n_{1}\left(\frac{k(k + 1)}{2}\right)$ as we wanted to show.
By assumption this quantity is non-negative, since it is the dimension of $\mathbf{L}_{d}(k^{n_{1}}).$

Finally, if we assume that $$\dim\mathbf{L}_{d}((k + 1)^{n_{1}}) = \dim\mathbf{L}_{k - 1}(m^{n_{2}}) = 0,$$ then $$\dim\mathbf{L}_d(k^{n_1})-\dim\mathbf{L}_d((k+1)^{n_1})+n_1(\dim\mathbf{L}_k(m^{n_2})-\dim\mathbf{L}_{k-1}(m^{n_2})) - n_{1}(k + 1)$$ 
is equal to 
$$\frac{d(d + 3)}{2} + 1 - n_{1}\frac{k(k + 1)}{2} + n_{1}\frac{k(k + 3)}{2} + n_{1} - n_{1}n_{2}\frac{m(m + 1)}{2} - n_{1}(k + 1).$$
We get precisely $\frac{d(d + 3)}{2} + 1 - n_{1}n_{2}\frac{m(m + 1)}{2}.$ 
\hfill $\Box$\\

We then use the above computations to find $\dim\mathbf{L}_{0}$ in several cases, in particular applying $(*)$ and the two claims.

\begin{Claim}
The above calculations imply
$$\dim\mathbf{L}_0=\frac{d(d+3)}{2}+1-n_1n_2\frac{m(m+1)}{2}$$
if one of the following conditions holds
$$\begin{array}{ll}
\hbox{i. } & \frac{d(d+3)}{2}+1-n_1\frac{(k+1)(k+2)}{2}\geq 0\hbox{ and }\frac{(k-1)(k+2)}{2}+1-n_2\frac{m(m+1)}{2}\geq 0\\
 & \\
\hbox{ii. } & \frac{d(d+3)}{2}+1-n_1\frac{(k+1)(k+2)}{2}\geq 0,\ \frac{(k-1)(k+2)}{2}+1-n_2\frac{m(m+1)}{2}\leq 0\\
 & \hbox{and  } \frac{k(k+3)}{2}+1-n_2\frac{m(m+1)}{2}\geq 0\\
 & \\
\hbox{iii. } & \frac{d(d+3)}{2}+1-n_1\frac{(k+1)(k+2)}{2}\leq 0,\ \frac{(k-1)(k+2)}{2}+1-n_2\frac{m(m+1)}{2}\geq 0\\
 & \hbox{and  } \frac{d(d+3)}{2}+1-n_1\frac{k(k+1)}{2}\geq 0\\
 & \\
\hbox{iv. } & \frac{d(d+3)}{2}+1-n_1\frac{(k+1)(k+2)}{2}\leq 0,\ \frac{(k-1)(k+2)}{2}+1-n_2\frac{m(m+1)}{2}\leq 0\\
 & \hbox{and  } \frac{d(d+3)}{2}+1-n_1n_2\frac{m(m+1)}{2}\geq 0
\end{array}$$
The same way we notice that
$\dim\mathbf{L}_0=0$ if all the inequalities
$$\begin{array}{cccc}
\hbox{a)} & \frac{d(d+3)}{2}+1-n_1\frac{(k+1)(k+2)}{2}\leq 0, & \\
 & \\ 
\hbox{b)} & \frac{(k-1)(k+2)}{2}+1-n_2\frac{m(m+1)}{2}\leq 0, &  \\
 & \\
\hbox{c)} & \frac{d(d+3)}{2}+1-n_1n_2\frac{m(m+1)}{2}\leq 0 & \hbox{ hold. } 
\end{array}$$
\end{Claim}
\prf
Suppose $\dim\mathbf{L}_{d}((k + 1)^{n_{1}})\geq 0$ and $\dim\mathbf{L}_{k - 1}(m^{n_{2}})\geq 0$ then 
$\dim\ker(\rho_Y)+\dim\ker(r_1,\ldots,r_{n_1})+\dim(\im(\rho_Y)\cap\im(r_1,\ldots,r_{n_1}))$
is equal to $$\frac{d(d + 3)}{2} + 1 - n_{1}n_{2}\frac{m(m + 1)}{2} - n_{1}(k + 1) + n_{1}(k + 1)$$ 
which is $\frac{d(d + 3)}{2} + 1 - n_{1}n_{2}\frac{m(m + 1)}{2}.$
We also note that, under the given assumptions, this quantity is non-negative since it is given as a sum of two dimensions.

If $\dim \mathbf{L}_{d}((k + 1)^{n_{1}})\geq 0, \dim\mathbf{L}_{k}(m^{n_{2}})\geq 0$ and $\dim \mathbf{L}_{k - 1}(m^{n_{2}}) = 0,$ then the sum
$$\frac{d(d + 3)}{2} + 1 - n_{1}\frac{(k + 1)(k + 2)}{2} + n_{1}\left(\frac{k(k + 3)}{2} + 1 - n_{2}
\frac{m(m + 1)}{2}\right)$$ 
gives $\frac{d(d + 3)}{2} + 1 - n_{1}n_{2}\frac{m(m + 1)}{2}$ as stated.

In the third case, namely $\dim\mathbf{L}_{d}(k^{n_{1}})\geq 0, \dim\mathbf{L}_{k - 1}(m^{n_{2}})\geq 0$ and $\dim \mathbf{L}_{d}((k + 1)^{n_{1}}) = 0,$ then the sum 
$\dim\ker(\rho_Y)+\dim\ker(r_1,\ldots,r_{n_1})+\dim(\im(\rho_Y)\cap\im(r_1,\ldots,r_{n_1}))$
is equal to 
$$\frac{d(d + 3)}{2} + 1 - n_{1}\frac{k(k + 1)}{2} + n_{1}\left(\frac{(k - 1)(k + 2)}{2} + 1 - n_{2}
\frac{m(m + 1)}{2}\right)$$ 
which again gives $\frac{d(d + 3)}{2} + 1 - n_{1}n_{2}\frac{m(m + 1)}{2}.$

If $\dim\mathbf{L}_{d}((k + 1)^{n_{1}}) = \dim\mathbf{L}_{k - 1}(m^{n_{2}}) = 0$ and also $$\frac{d(d + 3)}{2} + 1 - n_{1}n_{2}\frac{m(m + 1)}{2}\geq 0,$$ then, since 
$\dim\ker(\rho_Y)+\dim\ker(r_1,\ldots,r_{n_1})=0$
$$\begin{array}{c}
\dim\ker(\rho_Y)+\dim\ker(r_1,\ldots,r_{n_1})+\dim(\im(\rho_Y)\cap\im(r_1,\ldots,r_{n_1}))=\\
\\
\dim(\im(\rho_Y)\cap\im(r_1,\ldots,r_{n_1}))
= \frac{d(d + 3)}{2} + 1 - n_{1}n_{2}\frac{m(m + 1)}{2}. 
\end{array}$$

Finally, to prove the last statement,  both 
$\dim\ker(\rho_Y)+\dim\ker(r_1,\ldots,r_{n_1})$ and $\dim(\im(\rho_Y)\cap\im(r_1,\ldots,r_{n_1}))$
vanish so their sum is zero as well, proving that $\dim \mathbf{L}_{0} = 0.$\hfill $\Box$

\section{Proof of Theorem~\ref{thm2}}
\label{proof}
For the proof of Theorem~\ref{thm} fix numbers $d,\ m,\ n_1$ and $n_2$.\\

First assume that the virtual dimension 
$$\frac{d(d+3)}{2}+1-n_1n_2\frac{m(m+1)}{2}$$
of the vector space $\mathbf{L}_d(m^{n_1n_2})$ is positive. We will show that it equals to the dimension of  
$\mathbf{L}_0$. In Section~\ref{comp} we proved that the equality holds, if we can choose  $k$ that satisfies one of the conditions i), ii), iii) or iv). 

If  $\frac{d(d+3)}{2}+1-n_1\frac{(k + 1)(k + 2)}{2} \geq 0,$ then the only condition we need to check is that  $\frac{k(k+3)}{2}+1-n_2\frac{m(m+1)}{2}\geq 0.$ Then either i) or ii) will be satisfied.

In the other case  $\frac{d(d+3)}{2}+1-n_1\frac{(k + 1)(k + 2)}{2} \leq 0,$ then, since we're assuming that 
$\frac{d(d+3)}{2}+1-n_1n_2\frac{m(m+1)}{2} \geq 0,$ we just need to show that  $\frac{d(d+3)}{2}+1-n_1\frac{k(k + 1)}{2} \geq 0.$ Then either iii) or iv) will be true.\\

Define $k$ to be the integer between  
$$k_l=\frac{1}{2}\left(-3+\sqrt{1+4n_2m(m+1)}\right)\hbox{ and } 
k_u=\frac{1}{2}\left(-1+\sqrt{1+4n_2m(m+1)}\right).$$
Such an integer exists since $k_u=k_l+1.$ Then 
$\frac{k_l(k_l+3)}{2}+1-n_2\frac{m(m+1)}{2}=0$ together with $k_l\leq k$ implies
$$\frac{k(k+3)}{2}+1-n_2\frac{m(m+1)}{2}\geq 0.$$ 
The same way
$\frac{d(d+3)}{2}+1-n_1\frac{k_u(k_u+1)}{2}=\frac{d(d+3)}{2}+1-n_1n_2\frac{m(m+1)}{2}\geq 0$ together with
$k\leq k_u$ implies
$$\frac{d(d+3)}{2}+1-n_1\frac{k(k+1)}{2}\geq 0.$$
This proves that the chosen $k$ satisfies at least one set of inequalities i), ii), iii) or iv). \\

Next assume that the virtual dimension 
$$\frac{d(d+3)}{2}+1-n_1n_2\frac{m(m+1)}{2}$$
of the vector space $\mathbf{L}_d(m^{n_1n_2})$ is negative. We will show that in this case $\dim\mathbf{L}_0=0.$ 
Recall from the previous section that $\dim\mathbf{L}_0=0$ if we can find an integer $k$ such that all inequalities 
a), b), c)  hold. As before, define $k$ to be an integer between 
$$k_l=\frac{1}{2}\left(-3+\sqrt{1+4n_2m(m+1)}\right)\hbox{ and } 
k_u=\frac{1}{2}\left(-1+\sqrt{1+4n_2m(m+1)}\right).$$
The last inequality $c)$ is automatically fulfilled.
From
$$\frac{d(d+3)}{2}+1-n_1\frac{(k_l+1)(k_l+2)}{2}=\frac{d(d+3)}{2}+1-n_1n_2\frac{m(m+1)}{2}\leq 0$$
together with $k_l\leq k$ we get
$$\frac{d(d+3)}{2}+1-n_1\frac{(k+1)(k+2)}{2}\leq 0.$$
And $\frac{(k_u-1)(k_u+2)}{2}+1-n_2\frac{m(m+1)}{2}=0$ together with $k\leq k_u$ implies
$$\frac{(k-1)(k+2)}{2}+1-n_2\frac{m(m+1)}{2}\leq 0.$$\\

Therefore we proved, that we can always choose integer $k$ in the linear system $\mathbf{L}_0=\mathbf{L}_0(d,k,m,n_1n_2)$
such that 
$$\dim\mathbf{L}_0=\max\left\{\frac{d(d+3)}{2}+1-n_1n_2\frac{m(m+1)}{2},\ 0\right\}.$$ 
This finishes the proof of Theorem~\ref{thm2}.\\

\section{Applications}
\label{appl}
\subsection{A proof of Evain's theorem}

It is well known that the system of plane curves of degree $d$ passing through $4$ general points with homogeneous multiplicity $m$ is non-special. A proof can be found in ~\cite{cil_mir}. Given  this, Theorem ~\ref{thm}  implies  Evain's theorem ~\cite{evain}: namely that all systems of plane curves of degree $d$ through 
$4^{h}$ points with homogeneous multiplicity $m$ are non-special.

\subsection{The problem of $9^{h}$ points}

With similar methods, we can prove another theorem.
\begin{theorem}
The linear systems $\bl_d(m^{9^{h}})$ are non special.
\end{theorem}
As in the previous case, this theorem follows as a corollary from our main Theorem~\ref{thm}; since the case of $9$ points is known and completely analogous to the case of $4$ points.

\subsection{The problem of $4^{h}9^{k}$ points}

The following theorem is also a consequence of our Theorem~\ref{thm}.
\begin{theorem}
The linear systems $\bl_d(m^{4^{h}9^{k}})$ are non special for all integers $d, m, h$ and $k.$
\end{theorem}

\noindent {\small \sc Anita Buckley \\
\noindent Mathematics Institute, University of Warwick\\
\noindent Coventry CV4 7AL, United Kingdom \\
\noindent E-mail address: \tt mocnik@maths.warwick.ac.uk}\\
\noindent and\\
\noindent {\small \sc Marina Zompatori\\
\noindent Department of Mathematics, Boston University\\
\noindent Boston MA 02215, United States\\
\noindent E-mail address: \tt marinaz@math.bu.edu}\\

\begin{thebibliography}{99} 
{\small
\bibitem{cil_mir} {C. Ciliberto} and {R. Miranda}, {\it Degenerations of Planar Linear Systems}, J. Reine Angew. Math. 
{\bf 501} (1998), 191--220.
\bibitem{mir} {R. Miranda}, {\it Linear Systems of Plane Curves}, Notices of the AMS (2) {\bf 46} (1999), 192--202.
\bibitem{a_m} {A. Buckley} and {M. Zompatori}, {\it Generalization of the Transversality of the Restricted Systems}, to appear in Le Matematiche 
\bibitem{evain} {L. Evain}, {\it La Fonction de Hilbert de la r\'{e}union de $4^h$ points g\'{e}n\'{e}riques de $\bp^2$ de m\^{e}me multiplicit\'{e}}, J. of Alg. Geom. {\bf 8} (1999).
\bibitem{hir1} {A. Hirschowitz}, {\it Existence de faisceaux r\'{e}flexive de rang deux sur $\bp^3$ a bonne cohomologie},
Inst. Hautes \'{E}t. Sci. {\bf 66} (1988), 104--137.
\bibitem{hir2} {A. Hirschowitz}, {\it Une Conjecture pour la Cohomologie des Diviseurs sur les Surfaces Rationelles G\'{e}n\'{e}riques}, J. Reine Angew. Math. {\bf 397} (1989), 208--213.
}\end{thebibliography}
\end{document}